\documentclass[12pt]{article}
\usepackage[utf8]{inputenc}
\usepackage[polish,english]{babel}
\usepackage[T1]{fontenc}
\usepackage{indentfirst}
\usepackage{graphicx}
\usepackage{amsmath}
\usepackage{url}
\usepackage[breaklinks]{hyperref}
\usepackage{bbm} 
\newcommand{\pd}[2]{\frac{\partial #1}{\partial #2}}
\newcommand{\sfrac}[2]{{\scriptstyle \frac{#1}{#2}}}

\begin{document}

\title{A mathematical model of the Mafia game}
\author{Piotr Migda{\l}
\thanks{Institute of Theoretical Physics, University of Warsaw, Ho\.{z}a 69, 00-681 Warsaw, Poland}
\thanks{Faculty of Mathematics, Informatics and Mechanics, University of Warsaw, Banacha 2, 02-097 Warsaw, Poland}
\thanks{Current address: ICFO--Institut de Ci\`{e}ncies Fot\`{o}niques, Av. Carl Friedrich Gauss, num. 3, 08860 Castelldefels (Barcelona), Spain}
\thanks{E-mail address: pmigdal@gmail.com}
}

\maketitle

\begin{abstract}
Mafia (also called Werewolf) is a party game. The participants are divided into two competing groups: citizens and a mafia. The objective is to eliminate the opponent group. The game consists of two consecutive phases (day and night) and a certain set of actions (e.g. lynching during day). The mafia members have additional powers (knowing each other, killing during night) whereas the citizens are more numerous.

We propose a simple mathematical model of the game, which is essentially a pure death process with discrete time.  We find the closed-form solutions for the mafia winning-chance $w(n,m)$ as well as for the evolution of the game. Moreover, we investigate the discrete properties of results, as well as their continuous-time approximations.

It turns out that a relatively small number of the mafia members, i.e. proportional to the square root of the total number of players, gives equal winning-chance for both groups. Furthermore, the game strongly depends on the parity of the total number of players.
\end{abstract}

\section{Introduction}
\label{s:introduction}

Mafia (also called Werewolf) is a popular party game \cite{WikipediaMafia2010}. The participants collect in a circle and a game coordinator assigns each player to one of two groups: a mafia or citizens. Citizens know only their own identity, whereas mafia members know identity of their fellows. The objective is to eliminate the opponent group. The game consists of two alternating phases (day and night). During the day, all players can discuss and vote who they want to lynch. During the night the mafia kills a citizen of their choice.

There are many variants of the Mafia game. A typical modification of the gameplay is the addition of  characters with special abilities. For example, two common special citizens are: Detective (who checks every night if a chosen person is in the mafia) and Nurse (who may protect a victim from being killed, if she chooses correctly).

In this paper we analyze the simplest version of the Mafia game: with only mafia and citizens, and without any special players or additional rules. Even if it may not be the most popular variant, it is the best one for mathematical modeling. Despite the game having a complex psychological component it is possible to create a stochastic model of the Mafia game. That is, we focus on a game where all killings are at random, with decisions unaffected by the previous game course. It is a great simplification, and is going to follow only when participants fail to observe others' behavior, or make use of this information. However, even such a baseline model gives non-trivial predictions.

We consider a game beginning with $n$ players, out of which there are $m$ mafia members. The main questions we address are: 
\begin{itemize}
\item What is the probability that the mafia wins $w(n,m)$?
\item How do we expect the game to play out? That is, what is the chance that after a given time there is exactly a certain number of mafia members?
\end{itemize}
Besides the direct answers (i.e. closed form expressions), we study approximations, qualitative behavior and some special cases.

It turns out that a relatively small number of the mafia members , i.e. proportional to the square root of the total number of players, gives equal winning-chance for both groups. Furthermore, the game strongly depends on the parity of the total number of players.

There are only a few previous papers on the Mafia game. Work of Braverman, Etesami and Mossel \cite{Etesami2008} basically proposes the same model as the one presented in this paper. They calculate a simple asymptotic formula for the mafia-winning chance, $w(n,m)\propto m/\sqrt{n}$. We go a few steps further, finding the closed form expressions and parity-dependent asymptotic approximations and also analyzing the dynamical properties.

Other research projects concentrate on the psychological aspect of the Mafia game, in particular --- deceiving. Topics analyzed include: patterns of a deceiver's interruptions and voice parameters \cite{Chittaranjan2010}, movement of face and hands \cite{Xia2007} and usage of language \cite{Zhou2008}. Also, the Mafia game is applied as a test for some videoconferencing setups \cite{Batcheller2007}.


The research is motivated by two goals, besides the sheer fun of calculation and writing.  First, a good theoretical model may be useful for the investigation of the psychological aspect of the Mafia game. Experimental deviations from the idealized behavior may give a valuable insight into psychology of strategy choosing, manipulation, deceiving, following others and hiding identity. Second, when a small, but well-informed and powerful group fights against the majority, the nature of the process may be similar to that of the Mafia game --- e.g. in high-stake corruption, terrorism and illegal oppositions.

The paper is organized as follows. In Section~\ref{s:model} we introduce the rules of the Mafia game and propose a simple mathematical model. Section~\ref{s:one} contains a special case of play with only a single mafia member. The results are not only simple and didactic, but also useful for a more general case. Section~\ref{s:qualitative} gives simple results on the qualitative behavior of the mafia-winning chance. Section~\ref{s:evolution} presents the dynamical aspect of the Mafia game and gives closed-form expressions for dynamics of the game.  We consider both the discrete-time model and its continuous time approximation. In Section~\ref{s:general} we provide the exact result for the mafia-winning chance, together with its asymptotic approximation. Section~\ref{s:conclusion} concludes the work and gives insight into possible extensions and applications. Appendices \ref{appendixa} and \ref{appendixb} contain derivations of formulas found in Sec.~\ref{s:evolution}.

\section{Model}
\label{s:model}

Before proceeding to mathematics, we need to write down our arbitrarily chosen rules of  the Mafia game. Some assumptions are made to simplify the model, others purely for convenience. Note that this Section is only a hand-waving transition from the real-life game to its mathematical model.

\begin{itemize}
    \item The game needs $n$ players and one more person to coordinate it.
    \item At the beginning players are randomly divided into ($n-m$) citizens and $m$ mafia members.
    \item Mafia members know the identity of each other, citizens --- only their own.
    \item There are two alternating phases, day and night, which together comprise a turn.
    \item During day, there are two consecutive subphases:
    \begin{itemize}
        \item Debate. Everyone still alive can say anything related to accusing or defending.
        \item Voting. Everyone has one vote for who should be lynched. The player who gets the highest number of votes is eliminated (in case of a tie, a random 'winner' is eliminated). The victim's faction is revealed.
    \end{itemize}
    \item During the night:
    \begin{itemize}
        \item Mafia jointly decides who they want to kill (eliminate).
    \end{itemize}
    \item The game continues until there is only one group (either the citizens or the mafia) left. That group wins.
\end{itemize}

Stated briefly, during the day one player is eliminated (either a citizen or a mafia member) whereas during the night the mafia kills one citizen. Consequently, during a single turn the possible transitions are: $(n,m)\rightarrow(n-2,m)$ and $(n,m)\rightarrow(n-2,m-1)$. Their probabilities, in principle, depend on $n$, $m$ and the course of the play (i.e. previous debates and votings).
To start with, we need to get rid of the psychological aspect, restricting ourselves only to strategic and probabilistic parts of the Mafia game. It is a coarse approximation, as psychology plays an essential role in the Mafia game. Nevertheless, even such a bare model has interesting properties. Moreover, we want to neglect the debate phase as it is extremely difficult to formalize in a meaningful and useful way.

Now, let's make a hand-waving argument why the transition probability depends only on the current state $(n,m)$. Consider the scenario, in which players are lynched at random. The lynching victim is a member of the mafia with probability $m/n$, so
\begin{align}
P\left[ (n-2,m)| (n,m) \right] = \frac{n-m}{n} \quad \hbox{and} \quad P\left[ (n-2,m-1)| (n,m) \right] = \frac{m}{n}.\label{eq:transition}
\end{align}
Let's denote by $w(n,m)$ the mafia-winning chance with random lynchings, induced by the above transition probabilities.

Suppose that there is a better strategy (in the common sense of the word) for citizens, giving it $w'(n,m)< w(n,m)$. Then the mafia can force a random lynch, just by pretending they are citizens, that is, acting during the debate and the voting phases in the same way as citizens do. It is possible as:
\begin{itemize}
    \item Each mafia member knows more than any citizen.
    \item In such case there is no way to detect a mafia member.
\end{itemize}
Bear in mind that the second point works only in a game with no psychological part, when neither voting scheme nor behavior gives no trace of one's faction. For example, no-one is going to turn red when accused of being in the mafia.

Suppose that there is a better strategy for the mafia, $w''(n,m) > w(n,m)$. If the citizens are in majority, they can force a random lynching by the following procedure:
\begin{itemize}
    \item One citizen says "let every of us give a number $k_i$ from $1$ to $n$, we will kill the $\sum k_i \bmod n$-th".
    \item The number $\sum k_i \bmod n$ is random, so is the chosen player.  
    \item Every citizen to vote for the unlucky person whose number comes up.
    \item As there are more citizens than mafia members, the mafia votes do not play any role.
\end{itemize}

If the citizens are in the minority (i.e. $m>n-m$) no matter what the moves are, the mafia will win. For $m=n-m$ let's assume a random player is lynched.  If using explicit randomness is not forbidden --- instead of the sum-mod-$n$ game (similar to a child's counting-out game, when numbers are shown with fingers) a citizen can roll a $n$-sided die or provide $n$ straws (with one shorter, indicating the victim).
 
Consequently, random lynch is a well-justified strategy. A more detailed argument is shown in \cite{Etesami2008}.
Of course in the real world it would be really boring to play the Mafia game in such way. But well, since we dropped the psychological side, we neglect 'boredom' as well. Moreover, even if players do not play lynching at random, the effective choices may result in probabilities similar to that of our idealized model.
As a direct consequence of \eqref{eq:transition}, the mafia winning-chance can be expressed as a recurrence equation 
\begin{align}
    w(n,m)=\left\{
    \begin{array}{ll}
        0 & \hbox{if }m=0,\\
        1 & \hbox{if }m>n-m,\\
        \frac{n-m}{n}w(n-2,m)+\frac{m}{n}w(n-2,m-1)&\hbox{in all other cases}.\\
    \end{array}
    \right.\label{eq:w}
\end{align}
The above equation is the core of our paper. We investigate it for one mafia member, its qualitative behavior and closed-form solution along with its asymptotic formula.

It is important to point out that we may set different boundary conditions, i.e. in which if there is the same number of citizens and mafia members, mafia wins. That is, in \eqref{eq:w} we may write $(1 \hbox{ if }m\geq n-m)$ instead of $(1 \hbox{ if }m>n-m)$. In that case all results can be reproduced, perhaps
in a slightly changed form (i.e. different constants and even$\leftrightarrow$odd).  

\begin{figure}[!htbp]
    \centering
        \includegraphics[width=0.50\textwidth]{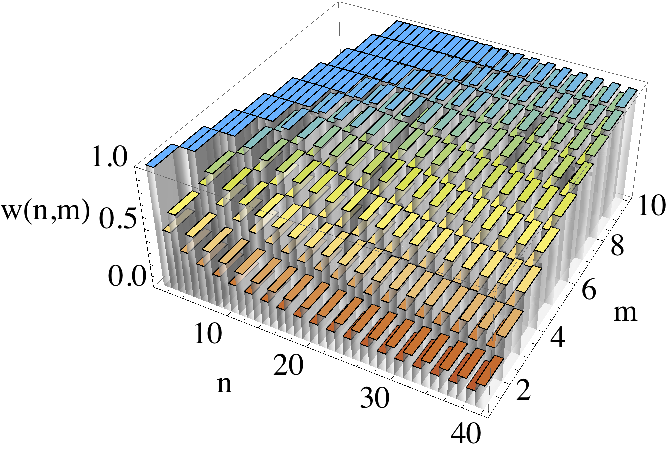}
    \caption{A numerical plot of the probability that the mafia wins $w(n,m)$. }
    \label{fig:wnm}
\end{figure}

\section{Game with one mafia member}
\label{s:one}

To start with, we want to consider a very simple case --- the game with only a single mafia member. That is, in this Section we analyze $w(n,1)$. Besides the simplicity, there are two more motivations: this case gives insight into some general properties and the results will be useful in the further part of the paper.

In a game with one mafia member, their winning probability is the chance that during every lynch a citizen is killed. Bearing in mind that we start the game with the day phase, we get
 \begin{align}
w(n,1)=\frac{n-1}{n}\cdot\frac{n-3}{n-2}\cdot \ldots \cdot \frac{1+(n \bmod 2)}{2+(n \bmod 2)}=\frac{(n-1)!!}{n!!},\label{eq:wn1}
\end{align}
where $n!!$ is the double factorial of $n$.
The above formula has an explicit dependence on the parity of $n$. Even though it is obvious that addition of $2$ citizens reduces mafia-winning chance $w(n+2,1)<w(n,1)$, it is not the case for addition of a single citizen. In fact, adding another player to make the total odd (holding $m=1$) increases the mafia-winning chance. To understand this let's use an example. In a play with one citizen and one mafia member, there is tie in voting, so lynching relies on tossing a coin and $w(2,1)=1/2$. In a game with two citizens, the mafia also wins after lynching a citizen (the second one is to be killed during the night). This time it is more difficult to hunt mafia, and $w(3,1)=2/3$.
One may think that we encountered an issue of boundary conditions (i.e. killing a random player in a tie) or that we considered too small $n$. Neither is the case. A careful reader may easily check that changing $(1$ if $m>n-m)$ to $(1$ if $m\geq n-m)$ in \eqref{eq:w} we get an analogous phenomena.

We may simply check that adding an odd player always increases the mafia winning-chance. Let's prove it by induction. As the base we use $w(3,1) > w(2,1)$. Then for the inductive assumption $w(2k+1, 1) > w(2k, 1)$ we need to show that $w(2(k+1)+1, 1)>w(2(k+1), 1)$:
\begin{align}
 w(2k+3, 1) &= \frac{2k+2}{2k+3} w(2k+1, 1) > \frac{2k+1}{2k+2} w(2k+1, 1)\label{eq:zabkidowod}\\
 &> \frac{2k+1}{2k+2} w(2k, 1) = w(2k+2, 1).\nonumber
\end{align}
Consequently, we have proved by mathematical induction that $w(2k+1,1)>w(2k,1)$ for every $k>0$.


Let's consider $w(n,1)$ averaged  over neighboring numbers. We use the geometric mean, as it simplifies the result
\begin{align}
\sqrt{w(n-1,1)w(n,1)}=\sqrt{\frac{(n-2)!!}{(n-1)!!}\frac{(n-1)!!}{n!!}}=\frac{1}{\sqrt{n}}.
\end{align}
The averaging may be considered an approximation in which we get a monotonic function (without the sawtooth pattern). To get some information on the dependence of $w(n,1)$ on the parity of $n$ we consider the following ratio 
\begin{align}
     \frac{w(2k+1,1)}{w(2k,1)}  &= \frac{(2k)!!}{(2k+1)!!} / \frac{(2k-1)!!}{(2k)!!} = 
 \left( \frac{(2k)!!}{(2k-1)!!} \right)^2\frac{1}{2k+1}\label{eq:prewallis}.
\end{align}
The above expression has a limit when $k$ goes to infinity, which can be found with the aid of the Wallis formula:
\begin{align}
    \frac{\pi}{2}=\lim_{k\rightarrow\infty} \left( \frac{(2k)!!}{(2k-1)!!} \right)^2\frac{1}{2k+1}.\label{eq:wallis}
\end{align}

We may develop approximate formulas for the single mafia member winning-chance, which take into account parity of the number of players 
\begin{align}
 w(2k,1) &= \left(\frac{w(2k+1,1)}{w(2k,1)}\right)^{-\frac{1}{2}} \sqrt{w(2k,1)w(2k+1,1)}\label{eq:m1przyblizenia}\\
    &\approx \sqrt{\frac{2}{\pi}} \frac{1}{\sqrt{2k+1}},\nonumber\\
 w(2k+1,1) &= \left(\frac{w(2k+1,1)}{w(2k,1)}\right)^{\frac{1}{2}} \sqrt{w(2k,1)w(2k+1,1)}\nonumber\\
  & \approx \sqrt{\frac{\pi}{2}} \frac{1}{\sqrt{2k+1}}.\nonumber
\end{align} 
Or, bearing in mind $(n+1)/n\rightarrow 1$,
\begin{align}
 w(n,1) &\approx \left(\frac{\pi}{2}\right)^{(n \bmod 2)-1/2} \frac{1}{\sqrt{n}}.\label{eq:m1approxnicer}
\end{align} 
 
As we see, the recurrence equation produces results that might be counter-intuitive. Imagine there are four of us and we are playing with a single mafia member. His winning-chance is $w(4,1)=\frac{3}{8}=0.375$. Then we invite five more friends to play as additional citizens. Consequently, the mafia winning-chance rise to $w(9,1)=\frac{128}{315}\approx 0.406$ (contrary to a naive expectation).
 
\begin{figure}[!htbp]
    \centering
        \includegraphics[width=0.35\textwidth]{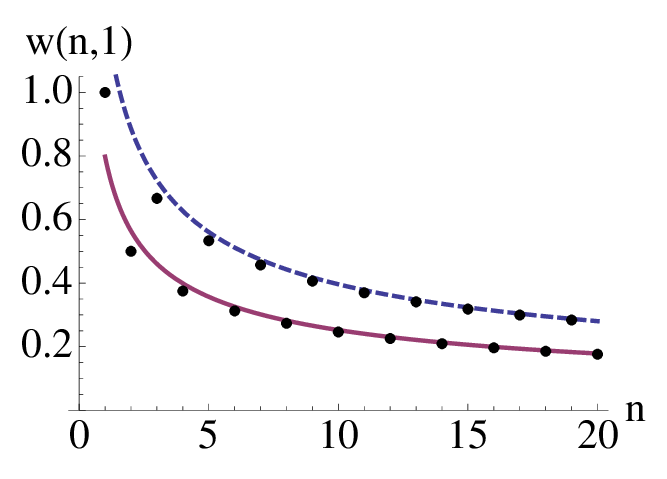}
    \caption{Plot of the mafia winning-chance for a game with a single mafia member, that is $w(n,1)$. Dots illustrate the exact result \eqref{eq:wn1}, whereas the lines represent the approximations \eqref{eq:m1approxnicer} --- solid and dashed line for even and odd number of players, respectively.}
    \label{fig:m1}
\end{figure}

\section{Qualitative properties of $w(n,m)$}
\label{s:qualitative}

Given the surprising effect of adding an additional citizen to a one-mafia-member game, some more questions naturally arise:
\begin{itemize}
    \item Does the addition of two citizens increase their winning-chance, i.e. is $w(n+2,m)<w(n,m)$?
    \item Does the addition of one citizen and one mafia member increase the mafia winning-chance, i.e. is $w(n+1,m+1)>w(n,m)$?
    \item Does the change of one player from a mafia member to a citizen decrease the mafia winning-chance, i.e. is $w(n,m+1)>w(n,m)$?
\end{itemize}
In fact, the above questions are equivalent, as we can see using using the recurrence relation \eqref{eq:w}
$$w(n,m)=\frac{n-m}{n}w(n-2,m)+\frac{m}{n}w(n-2,m-1)$$
for $n-m\geq m>0$. After some cosmetic arithmetical operations we get
\begin{align}
    (n-m)\left[ w(n-2,m)-w(n,m) \right] &= m \left[ w(n,m)-w(n-2,m-1) \right],\label{eq:dwie_sciezki}\\
    n\left[ w(n,m) - w(n-2,m-1)\right] &= (n-m) \left[ w(n-2,m)-w(n-2,m-1) \right].\nonumber
\end{align}
The above relations can be directly translated into
\begin{align}
  w(n-2,m)>w(n,m)
     \Leftrightarrow  w(n,m)>w(n-2,m-1)
     \Leftrightarrow  w(n-2,m)>w(n-2,m-1).\label{eq:rownowaznie}
\end{align}

In addition to the equivalence of the 'obvious' inequalities, we are interested whether they are fulfilled. Let's use the third part of \eqref{eq:rownowaznie} as a conjecture, formulating in the following way:

If a mafia member is changed into a citizen, the mafia winning-chance does not increase, that is $w(n+1,m+1)\geq w(n+1,m)$. If additionally $n-m\geq m$, the inequality is strong, $w(n+1,m+1)> w(n+1,m)$.

Let's introduce an agent to a game with $n$ players among which there are $m$ mafia members. The agent is an additional player who simulates either a citizen or a mafia member. We are going to show that the agent does not have to reveal its identity until the very late part of the game.

As citizens are indistinguishable, let's assume that mafia kills them during night in a fixed sequence. Let the agent have the last position. During the lynching there are three possibilities 
\begin{itemize}
    \item A citizen dies.
    \item A mafia member dies.
    \item The agent dies.
\end{itemize}
In the case of the agent's death, his affiliation has no meaning. In the two remaining  situations, the game continues with the agent. The agent's affiliation plays a decisive role in two situations:
\begin{itemize}
    \item When all mafia members are dead.
    \item When $0$ or $1$ citizens and $1$ or $2$ mafia members remain.
\end{itemize}
Note that in a real game, the agent's previous voting may be important in creating or breaking ties, not to mention the debate phase. But random voting squashes all of these considerations. In both cases when the agent is a mafia member, the mafia winning-chance increase:
\begin{center}
\begin{tabular}{|c|c|c|}\hline
        Game state $(n,m)$ & $w$ if the agent is of citizens & $w$ if the agent is of mafia  \\ \hline \hline
        $(n,0)$ &   $0$                         &   $w(n,1)$    \\ \hline
        $(1,1)$ &   $\frac{1}{2}$       &   $1$             \\ \hline
        $(2,1)$ &   $\frac{1}{3}$       &   $1$             \\ \hline
        $(3,2)$ & $\frac{1}{4}$     &   $1$             \\ \hline
\end{tabular}
\end{center}

The total mafia winning-chance is equal to the mean winning-chance, averaged over disjoint games. The weight is the probability of achieving a such game. The change of the agent's affiliation does not affect the weights, but in the ending always increases the mafia winning-chance, $w(n+1,m+1)\geq w(n+1,m)$, thus proving our first conjecture. Additionally, when the citizens outnumber the mafia, there is a game in which all mafia members are killed and the agent's affiliation plays a role, hence  $w(n+1,m+1)> w(n+1,m)$, which was our second conjecture

Let's show one more property --- when $n-m\geq m$ adding an odd player always increases the mafia winning-chance, not only in games with one mafia member \eqref{eq:zabkidowod}. We will use mathematical induction, using the boundary conditions \eqref{eq:w} and the  already shown properties of $w(n,1)$ as the basis. The inductive assumption is that the property holds for the total number of players $n\leq 2k+1$, and the inductive step is the following:
\begin{align}
w(2(k+1)+1,m) &= \frac{2k+3-m}{2k+3}w(2k+1,m)+\frac{m}{2k+3}w(2k+1,m-1)\label{eq:paritygen}\\
 &> \frac{2k+3-m}{2k+3}w(2k,m)+\frac{m}{2k+3}w(2k,m-1)\nonumber\\
 &> \frac{2k+2-m}{2k+2}w(2k,m)+\frac{m}{2k+2}w(2k,m-1)\nonumber\\
 &= w(2(k+1),m),\nonumber
\end{align}
where we subsequently used: the inductive assumption, the previously proven relation $w(2k,m)>w(2k,m-1)$ and the fact that (when averaging two components) putting less weight on the greater term makes the sum lower.

\section{Evolution}
\label{s:evolution}

In two previous Sections we considered the mafia winning-chance with respect to the initial number of players ($N$) and mafia members ($M$). It may be interesting to explore dynamics of the play, that is, analyze the probability $p_m(t)$ that after the $t$-th turn there will be exactly $m$ mafia members. Even though this goes beyond the \eqref{eq:w}, it is a direct consequence of the random lynching model \eqref{eq:transition}. Throughout this section we use capital letters $N$ and $M$ to denote the initial conditions.

Note that every turn, or day and night, two players are killed, so the total number of players decreases $n(t)=N-2t$.

\subsection{Discrete time}
\label{s:discrete}

The evolution of each $p_m(t)$ is governed by a set of recurrent equations
\begin{align}
    p_m(t+1)  = \frac{(N-2t)-m}{N-2t}p_m(t)+\frac{m+1}{N-2t}p_{m+1}(t),\label{eq:pm_dysk_ewol}
\end{align}  
where $m\in \mathbbm{N}$ and with the initial condition $p_m(0)=\delta_{mM}$.
 
Note that the above is a stochastic process with $M+1$ distinct states and $(m)\rightarrow(m-1)$. In general, such a Markov process is called a pure death process (a subclass of birth and death processes). In our case the transition probabilities change with time, which is a slight complication.

It is essential to point that the equations in \eqref{eq:pm_dysk_ewol} are correct only when $m\leq n(t)$, that is, when there is at least one citizen alive and the game can be still played. Otherwise they have no meaning. But why does an erroneous equation for  $m > n(t)$ not spoil states with $m\leq n(t)$?  For a given time $t$ let $p_m(t)$ be the erroneous probability with the smallest index, i.e. $m=N-2t+1$. In the next turn it affects $p_{m-1}(t+1)$. But, now the smallest wrong index is $N-2(t+1)+1=m-2<m-1$. So the error propagates more slowly than the exclusion of the states.

The closed form solution for $p_m(t)$ reads
\begin{align}
    p_m(t)=\sum_{i=m}^M {M \choose i}{i \choose m} (-1)^{i-m} \frac{(N-2t)!!}{N!!}\frac{(N-i)!!}{(N-2t-i)!!}\label{eq:pm_dysk_rozw}
\end{align}
and its derivation is presented in Appendix \ref{appendixa}. An example is plotted in Fig. \ref{fig:disccont} (a).

The expression above, although exact, is fairly complicated. Let's analyze the mean number of mafia players after $t$ days and nights, that is
\begin{align}
    \langle m \rangle(t)=\sum_{m=0}^M m p_m(t).\label{eq:srednio_definicja}
\end{align}
In birth and death processes with transition coefficients proportional to state's label (as radioactive decay process or Yule process) such mean may be obtained by a straightforward calculation with the use of \eqref{eq:pm_dysk_ewol}
\begin{align}
    \langle m \rangle(t) &= \frac{N-2t+1}{N-2t+2} \langle m \rangle(t-1) = \left(\prod_{i=0}^{t-1}\frac{N-2i-1}{N-2i}\right) M.\label{eq:sred}
\end{align}
The result is correct only when $N-2t-M\geq0$, that is, the average is taken over probabilities of the correct states.

\subsection{Continuous time approximation}
\label{s:continuous}

Usually differential equations are simpler than difference equations. In this Section we crudely approximate the discrete-time evolution  \eqref{eq:pm_dysk_ewol}  by its continuous-time version. We know that
\begin{align}
    p_m(t+1)-p_m(t)\xleftarrow[\Delta=1]{}
    \frac{p_m(t+\Delta)-p_m(t)}{\Delta}
    \xrightarrow[\Delta\rightarrow0]{}\frac{d}{dt}p_m(t).\label{eq:rrrr}
\end{align}
Let's change in \eqref{eq:pm_dysk_ewol} the difference $p_m(t+1)-p_m(t)$ into the differential $\frac{d}{dt}\tilde{p}_m(t)$ 
\begin{align}
    \frac{d}{dt}\tilde{p}_m(t) =
     -\frac{m}{N-2t}\tilde{p}_m(t)
     +\frac{m+1}{N-2t}\tilde{p}_{m+1}(t),\label{eq:pm_ciag_ewol}
\end{align}
where $m\in \mathbbm{N}$ and with the initial condition $\tilde{p}_m(0)=\delta_{mM}$. When a relative change of the function derivative is over unit length interval $[t,t+1]$, the approximation should be good $\tilde{p}_m(t)\approx p_m(t)$.  Intuitively speaking, when the chance of lynching a mafia member during one day is small, such an approximation is justified. However, in the course of this paper we do not estimate the approximation error. 

\begin{figure}[!htbp]
    \centering
        \begin{tabular}{cc}
        \includegraphics[width=0.50\textwidth]{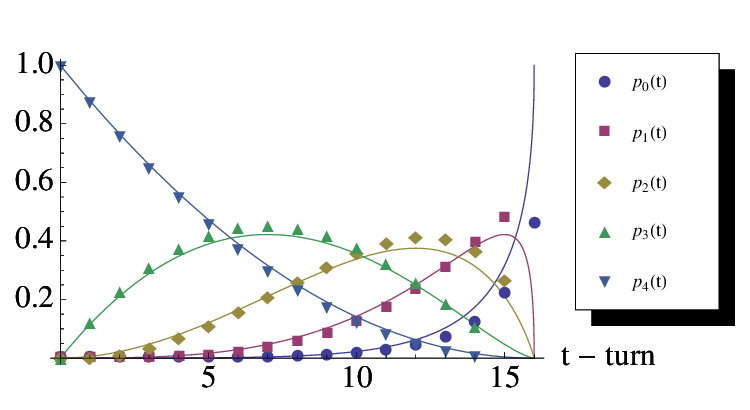} &
        \includegraphics[width=0.40\textwidth]{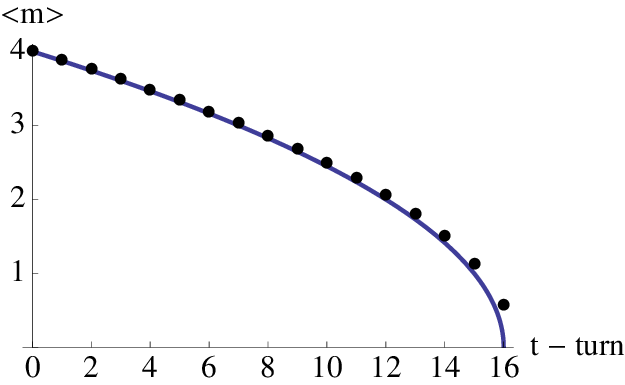}\\ 
        (a) & (b)
        \end{tabular}
    \caption{Dynamics of the Mafia game for the initial conditions $N=32$ and $M=4$. We compare results for the discrete-time evolution with its continuous-time counterparts. (a) The probability that after $t$ turns there are exactly $m$ mafia members alive --- the exact results $p_m(t)$ (dots) and their approximations $\tilde{p}_m(t)$ (lines). (b) The mean number mafia members --- the discrete result $\langle m \rangle(t)$ (dots) and its continuous time variant $\langle \tilde{m} \rangle(t)$ (lines).}
    \label{fig:disccont}
\end{figure}

The solution of the differential equation  \eqref{eq:pm_ciag_ewol} reads
\begin{align}
    \tilde{p}_m(t)={M \choose m}\left( 1-\sqrt{1-\frac{2t}{N}}\right)^{M-m}\left( \sqrt{1-\frac{2t}{N}} \right)^m\label{eq:pm_ciag_rozw},
\end{align}
see Appendix \ref{appendixb} for a detailed derivation. This is more convenient to work with than the exact expression \eqref{eq:pm_dysk_rozw}. An example is plotted in Fig. \ref{fig:disccont} (a).

For example, maxima of $\tilde{p}_m(t)$ are easily found
\begin{align}
    t_m = \frac{N}{2}\left(1-\left(\frac{m}{M}\right)^2 \right). \label{eq:pm_maksima}
\end{align}

The mean number of mafia players, defined as in \eqref{eq:srednio_definicja}, is
\begin{align}
    \langle \tilde{m} \rangle(t)= M \sqrt{1-\frac{2t}{N}}\label{eq:srednio_przybl},
\end{align}
which is a simpler expression than its analogue \eqref{eq:sred}. The comparison of the discrete-time formula with its continuous-time approximation is in Fig. \ref{fig:disccont} (b).

It is tempting to estimate the mafia winning-chance. As $\tilde{p}_0(t)$ is the probability that the citizens won, setting time for which there is only one player left gives the total citizen winning-chance:
\begin{align}
    \tilde{w}(n,m)&=1-\tilde{p}_0\left(\sfrac{n-1}{2}\right)=1-\left( 1-\sqrt{1-\frac{n-1}{n}}\right)^m\\
    &\approx \frac{m}{\sqrt{n}}.\nonumber
\end{align}
Not surprisingly, in the continuous-time approximation, there is no explicit dependence on parity (the discrete properties are lost).

\section{General solution}
\label{s:general}

The result for $p_m(t)$ \eqref{eq:pm_dysk_rozw} is not only interesting by itself --- it also gives us the solution for the mafia winning-chance \eqref{eq:w}. Note that $p_0(t)$ is the probability that after $t$ turns citizens win. When total number of players is odd, game lasts for $t_{odd}=(N-1)/2$ days and nights, leaving 1 player alive. Thus $p_0(t_{odd})$ is the probability the citizen win the game. For an even number of players, after $t_{even}=N/2$ turns there are no players alive. But now $p_0(t_{even})$ is the probability that there are no mafia members, so in fact --- no-one was killed during the last night and citizen win. So $w(n,m)=1-p_0(\sfrac{N-N\bmod 2}{2})$. Consequently, the closed form formula for the mafia winning-chance reads
\begin{align}
    w(n,m)&=1-\sum_{i=0}^m {m \choose i} (-1)^i \frac{(n-i)!!}{n!!((n\bmod 2)-i)!!}.\label{eq:w_rozwiazanie}
\end{align}

Let's show what we have proven for a game with only one mafia member --- the asymptotic behavior \eqref{eq:m1approxnicer}. 

To get the asymptotic formula for $w(n,m)$ it suffices to notice that in the limit $n\rightarrow\infty$ (taken over a selected parity) only the first two terms of the sum contribute. For $i=0$ the term is equal to $1$. The terms with $i>1$ do not matter as $(n-i)!!/(n-1)!!\rightarrow 0$. Consequently, we may write
\begin{align}
    w(n,m)&\approx m \frac{(n-1)!!}{n!!}\label{eq:genapprox}\\
    &\approx \left(\frac{\pi}{2}\right)^{(n \bmod 2)-1/2} \frac{m}{\sqrt{n}}\nonumber,
\end{align}
where we used the result from game with one mafia member, that is, approximating \eqref{eq:wn1} with \eqref{eq:m1approxnicer}. The approximate formula \eqref{eq:genapprox} can be formally expressed as $\frac{w(n,m)}{m/\sqrt{m}}$

\begin{figure}[!htbp]
    \centering
        \includegraphics[width=0.60\textwidth]{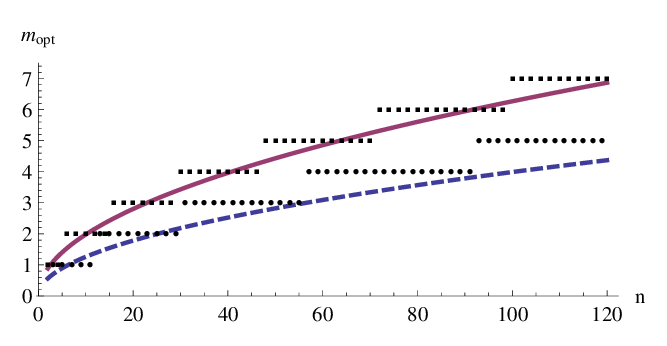}
    \caption{Optimal number of mafia members for a given number of players. Points show numerical results, that is, number of mafia members $m_{opt}$ for which $w(n,m_{opt})$ is the closest to $1/2$. Lines show the approximation \eqref{eq:optimalmafia}. Squares and solid line are for even number of players, whereas dots and dashed line --- odd.}
    \label{fig:mopt}
\end{figure}

It is tempting to ask for the optimal number of mafia members $m_{opt}$, for a given number of players. An interesting game is one in which both groups have the same chance of winning
\begin{align}
    m_{opt}& \approx \frac{1}{2}\left(\frac{\pi}{2}\right)^{-(n \bmod 2)+1/2} \sqrt{n}, \label{eq:optimalmafia}
\end{align}
which are plotted in Fig. \ref{fig:mopt}. Note that it is only a hand-waving argument, as we applied an asymptotic formula to a finite $n$, without investigating how does the convergence depend on $m$. Fortunately, for even $n$ the component with $i=2$ in the series \eqref{eq:w_rozwiazanie} vanishes as $1/(-2)!!=0$. Consequently, for $m^2$ of order of $n$ (as in \eqref{eq:optimalmafia}) the remaining part goes to zero as $n$ goes to infinity. Unfortunately, one cannot make the same argument for odd $n$ and for such case the approximation \eqref{eq:optimalmafia} may not hold.

\section{Further remarks}
\label{s:conclusion}

In every realistic Mafia game, there is a psychological part. A dry mathematical model may be a backbone of a more complex one, but certainly, is not enough to describe a real game. During course of play, citizens gain some information --- either by discovering another's identity by themselves, or trying to catch messages. Furthermore, voting may be subject to some kind of witch-hunt mentality. Moreover, rarely are all players the same --- usually there are ones with higher and lower influence on the others.

All these processes may be investigated mathematically, or numerically. There is still a lot of challenge in finding an appropriate model and obtaining the results. Comparison of the theory with the experiment may be crucial. 

Out of author's personal experience, in games when players don't know each other, usually a random person is killed (e.g. a person who is always nervous, not just because of being in the mafia). When they know each other very well, very dim signals can be used to reveal one's identity.

\subsection*{Acknowledgments}

The paper is based on my bachelor's thesis \cite{Migdal2009}, written under supervision of prof. Jacek Mi\k{e}kisz.

I would like to thank Michael Kleber, G. John Lapeyre, Marcin Kotowski and an anonymous reviewer for valuable remarks and language corrections. Additionally, I am grateful to all with whom I have played the Mafia game, with a special emphasis on friends from the Polish Children's Fund. Personally, I am a fan of Ktulu, a rather complex variant of the Mafia game.  In Ktulu there are 4 distinct factions, every player has a different special ability and the winning conditions rely on a special item.

\appendix

\section{Solution of difference equation for $p_m(t)$}\label{appendixa}

We want to solve the following set of  difference equations \eqref{eq:pm_dysk_ewol},
\begin{align}
    p_m(t+1)  = \frac{(N-2t)-m}{N-2t}p_m(t)+\frac{m+1}{N-2t}p_{m+1}(t)\quad : \quad m\in \mathbbm{Z} \nonumber
\end{align}
with the initial condition $p_m(0)=\delta_{mM}$. Let's introduce a generating function
\begin{align}
    F(t,z)=\sum_{m=0}^\infty p_m(t)z^m\label{eq:tworzaca_dysk},
\end{align}
with the respective initial condition $F(0,z)=z^M$. After side-by-side multiplication of \eqref{eq:pm_dysk_ewol} by $z^m$ and summation we get
\begin{align}
    F(t+1,z)&=F(t,z)-\frac{z}{N-2t}\pd{F(t,z)}{z}+\frac{1}{N-2t}\pd{F(t,z)}{z}\label{eq:tworz_lin}\\
    &=\underbrace{\left( 1-\frac{z-1}{N-2t}\pd{ }{z}\right)}_{\mathbf{B}(t)}F(t,z).\nonumber
\end{align}
The linear operator $\mathbf{B}(t)$ has the eigenvectors of the form $(z-1)^k$ with the corresponding eigenvalues $(1-k/(N-2t))$. As the eigenvectors do not depend on $t$, the computation of $F(t,z)$ is straightforward: 
\begin{align}
    F(t,z)&=\mathbf{B}(t-1)\mathbf{B}(t-2)\cdots\mathbf{B}(0) z^M\label{eq:tworz}\\
    &=\sum_{i=0}^M {M \choose i} (z-1)^i \left( 1-\frac{i}{N-2(t-1)}\right)\left( 1-\frac{i}{N-2(t-2)}\right)\cdots\left( 1-\frac{k}{N}\right)\nonumber\\
    &=\sum_{k=0}^M {M \choose i} (z-1)^i \frac{(N-2t)!!}{N!!}\frac{(N-i)!!}{(N-2t-i)!!}.\nonumber
\end{align}
Applying once more the binomial theorem, and comparing the result \eqref{eq:tworz} to the definition of the generating function \eqref{eq:tworzaca_dysk}  we get 
\begin{align}
    p_m(t)=\sum_{i=m}^M {M \choose i}{i \choose m} (-1)^{i-m} \frac{(N-2t)!!}{N!!}\frac{(N-i)!!}{(N-2t-i)!!},\nonumber
\end{align}
or \eqref{eq:pm_dysk_rozw}.

\section{Solution of differential equation for $\tilde{p}_m(t)$}\label{appendixb}

We want to solve the following set of differential equations \eqref{eq:pm_ciag_ewol},
\begin{align}
    \frac{d}{dt}\tilde{p}_m(t) =
     -\frac{m}{N-2t}\tilde{p}_m(t)
     +\frac{m+1}{N-2t}\tilde{p}_{m+1}(t)\quad : \quad m\in \mathbbm{Z},\nonumber
\end{align}
with the initial condition $\tilde{p}(0)=\delta_{mM}$. Once again, we employ a generating function 
\begin{align}
    G(t,z)=\sum_{m=0}^\infty \tilde{p}_m(t)z^m\label{eq:tworzaca_ciag},
\end{align}
with the respective initial condition $G(0,z)=z^M$. We get the partial differential equation
\begin{align}
    \pd{G(t,z)}{t}=\frac{(-z+1)}{N-2t}\pd{G(t,z)}{z},\label{eq:tworz_ciag_rown}
\end{align}
which we want to solve by the method of characteristics. Let's have integral curves of $G(t,z)$ in the form of $(t(\varphi),z(\varphi))$,
\begin{align}
    \pd{G(t,z)}{t}\frac{dt}{d\varphi}+\pd{G(t,z)}{z}\frac{dz}{d\varphi}=0\label{eq:char}.
\end{align}
Comparison of \eqref{eq:tworz_ciag_rown} with \eqref{eq:char} results in
\begin{align}
    \frac{dt}{N-2t}=d\varphi=\frac{dz}{z-1},
\end{align}
or
\begin{align}
    G(t,z)&=f\left(\ln\left|\sqrt{N-2t}(1-z)\right|\right).\label{eq:poziomicowo}
\end{align}
All we need is to find a function $f(x)$ that satisfies the initial condition. The explicit form of the generating function reads
\begin{align}
    G(t,z)=\left(1-\sqrt{1-\sfrac{2t}{N}}(1-z) \right)^M,
\end{align}
which gives the result
\begin{align}
    \tilde{p}_m(t)={M \choose m}\left( 1-\sqrt{1-\frac{2t}{N}}\right)^{M-m}\left( \sqrt{1-\frac{2t}{N}} \right)^m\nonumber
\end{align}
or \eqref{eq:pm_ciag_rozw}.

\end{document}